\documentclass[reqno,11pt]{amsart}
\usepackage{amsmath, latexsym, amsfonts, amssymb, amsthm, amscd}
\usepackage{graphics,epsf,psfrag}
\numberwithin{equation}{section}

\newtheorem{theorem}{Theorem}[section]
\newtheorem{lemma}[theorem]{Lemma}
\newtheorem{proposition}[theorem]{Proposition}

\newtheorem{rem}[theorem]{Remark}

\DeclareMathOperator{\sign}{\mathrm{sign}}

\newcommand{\ind}{\mathbf{1}}

\newcommand{\Z}{\mathbb{Z}}
\newcommand{\N}{\mathbb{N}}
\renewcommand{\tilde}{\widetilde}
\renewcommand{\hat}{\widehat}

\newcommand{\bP}{{\ensuremath{\mathbf P}} }
\newcommand{\bE}{{\ensuremath{\mathbf E}} }

\DeclareMathSymbol{\leqslant}{\mathalpha}{AMSa}{"36} 
\DeclareMathSymbol{\geqslant}{\mathalpha}{AMSa}{"3E} 
\DeclareMathSymbol{\eset}{\mathalpha}{AMSb}{"3F}     
\newcommand{\dd}{\,\text{\rm d}}             
\newcommand{\sumtwo}[2]{\sum_{\substack{#1 \\ #2}}} 
\newcommand{\supthree}[3]{\sup_{\substack{#1 \\ #2 \\ #3}}} 

\newcommand{\bbE}{{\ensuremath{\mathbb E}} }

\newcommand{\bbP}{{\ensuremath{\mathbb P}} }

\newcommand{\ga}{\alpha}

\newcommand{\gga}{\gamma}            

\newcommand{\gep}{\varepsilon}       

\newcommand{\go}{\omega}

\newcommand{\gl}{\lambda}

\makeatletter
\def\captionfont@{\footnotesize}
\def\captionheadfont@{\scshape}

\long\def\@makecaption#1#2{%
  \vspace{2mm}
  \setbox\@tempboxa\vbox{\color@setgroup
    \advance\hsize-6pc\noindent
    \captionfont@\captionheadfont@#1\@xp\@ifnotempty\@xp
        {\@cdr#2\@nil}{.\captionfont@\upshape\enspace#2}%
    \unskip\kern-6pc\par
    \global\setbox\@ne\lastbox\color@endgroup}%
  \ifhbox\@ne 
    \setbox\@ne\hbox{\unhbox\@ne\unskip\unskip\unpenalty\unkern}%
  \fi
  \ifdim\wd\@tempboxa=\z@ 
    \setbox\@ne\hbox to\columnwidth{\hss\kern-6pc\box\@ne\hss}%
  \else 
    \setbox\@ne\vbox{\unvbox\@tempboxa\parskip\z@skip
        \noindent\unhbox\@ne\advance\hsize-6pc\par}%
\fi
  \ifnum\@tempcnta<64 
    \addvspace\abovecaptionskip
    \moveright 3pc\box\@ne
  \else 
    \moveright 3pc\box\@ne
    \nobreak
    \vskip\belowcaptionskip
  \fi
\relax
}
\makeatother
\def\writefig#1 #2 #3 {\rlap{\kern #1 truecm
\raise #2 truecm \hbox{#3}}}

\newcommand{\tf}{\textsc{f}}

\begin{document}

\title[Coarse graining and critical slope]{Coarse graining,
  fractional moments and the critical slope of random copolymers}

 \author{Fabio~Lucio~Toninelli}
 \address{
CNRS and 
 ENS Lyon, Laboratoire de Physique,
 46 All\'ee d'Italie, 69364 Lyon, France\\ \\
\textit{E-mail address:} fabio-lucio.toninelli@ens-lyon.fr}

 \begin{abstract} For a much-studied model of random copolymer at a
   selective interface we prove that the slope of the critical curve
   in the weak-disorder limit is strictly smaller than $1$, which is
   the value given by the annealed inequality. The proof is based on a
   coarse-graining procedure, combined with upper bounds on the
   fractional moments of the partition function.  
\textit{Submitted to EJP on June 10, 2008,
final version accepted on January 30, 2009}
 \end{abstract}

\maketitle

\newpage
\section{Introduction}

We consider a model of {\sl copolymer at a selective interface}
introduced in \cite{cf:GHLO}, which has  attracted much
attention among both theoretical physicists and probabilists (we refer to
\cite{cf:GB} for general references and motivations).

Let $S:=\{S_n\}_{n\ge0}$ be the symmetric Simple Random Walk on $\Z$ started
at $S_0=0$, with law $\bP^{\rm SRW}$ such that the increments
$\{S_i-S_{i-1}\}_{i\ge0}$ are IID and $\bP^{\rm SRW}(S_1=\pm1)=1/2$.
The partition function of the model of size $N$ is given by
\begin{eqnarray}
  \label{eq:Zsrw}
  Z_{N,\go}:=\bE^{\rm SRW}\left[e^{-\gl\sum_{n=1}^N(\go_n+h)(1-\sign(S_n))}
\ind_{\{S_N=0\}}
\right]
\end{eqnarray}
where $h\ge0$, $\gl\ge0$ and $\{\go_n\}_{n\ge1}$ is a sequence of IID
standard Gaussian random variables (the quenched disorder).  We adopt
the convention that, if $S_n=0$, then $\sign(S_n):=\sign(S_{n-1})$.
One interprets $\gl$ as the inverse temperature (or coupling strength)
and $h$ as an ``asymmetry parameter'': if $h>0$, since the $\go_n$'s
are centered, the random walk overall prefers to be in the upper
half-plane ($S\ge0$).  It is known that the model undergoes a
delocalization transition: if the asymmetry parameter exceeds a
critical value $h_c(\gl)$ then the fraction of ``monomers'' $S_n$,
$n\le N$, which are in the upper half-plane tends to $1$ in the
thermodynamic limit $N\to\infty$ (delocalized phase), while if
$h<h_c(\gl)$ then a non-zero fraction of them is in the lower
half-plane (localized phase). What mostly attracts attention is the
slope, call it $m_c$, of the curve $\gl\mapsto h_c(\gl)$ in the limit
$\gl\searrow 0$: $m_c$  is expected to be a universal quantity, i.e.,
independent of the details of the law $\bP^{\rm SRW}$ and of the
disorder distribution (see next section for a more extended discussion
on this point). Already the fact that the limit slope is well-defined
and positive is highly non-trivial \cite{cf:BdH}.

Until now, all what was known rigorously about $m_c$ is that $2/3\le
m_c\le 1$, but numerically the true value seems to be rather around
$0.83$ \cite{cf:CGG}. The upper bound comes simply from annealing,
i.e., from Jensen's inequality, as explained in next section.  Our
main new result is that $m_c$ is strictly smaller than $1$.  The proof
works through a coarse-graining procedure in which one looks at the
system on the length-scale $k(\gl,h)$, given by the inverse of the
annealed free energy. The other essential ingredient is a
change-of-measure idea to estimate fractional moments of the partition
function (this idea was developed in \cite{cf:GLT} and \cite{cf:DGLT},
and used in the context of copolymers in \cite{cf:BGLT}).
Coarse-graining schemes, implemented in a way very different from
ours, have already played an important role in this and related
polymer models; we mention in particular \cite{cf:BdH}, \cite{cf:BCT}
and \cite{cf:AZy}.

\section{The general copolymer model}

As in \cite{cf:BGLT}, we consider a more general copolymer model which
includes \eqref{eq:Zsrw} as a particular case. Since the critical
slope is not proven to exist in this general setting, Theorem
\ref{th:slope} will involve a $\limsup$ instead of a limit.

Consider a renewal process $\tau:=\{\tau_0,\tau_1,\ldots\}$ of law
$\bP$, where $\tau_0:=0$ and $\{\tau_i-\tau_{i-1}\}_{i\in\N}$ is an IID
sequence of integer-valued random variables. We call
$K(n):=\bP(\tau_1=n)$ and we assume that $\sum_{n\in\N}K(n)=1$ (the
renewal is recurrent) and that $K(\cdot)$ has a power-law tail:
\begin{eqnarray}
\label{eq:K}
  K(n)\stackrel{n\to\infty}\sim \frac {{\mathcal C}_K}{n^{1+\ga}}
\end{eqnarray}
with $\ga>0$ and ${\mathcal C}_K>0$. As usual, the notation
$A(x)\stackrel{x\to x_0}\sim B(x)$ is understood to mean that $\lim_{x\to
  x_0}A(x)/B(x)=1$.

The copolymer model we are going to define depends on two parameters
 $\gl\ge0$ and $h\ge0$, and on a sequence
$\go:=\{\go_1,\go_2,\ldots\}$ of IID standard Gaussian random
variables (the quenched disorder), whose law is denoted by $\bbP$.
For a given system size $N\in\N$ and disorder realization $\go$, we
define the partition function
 $ Z_{N,\go}:=Z_{N,\go}(\gl,h)$ as
\begin{eqnarray}
  \label{eq:Z}
  Z_{N,\go}:=\bE\left[\prod_{j\ge 1:\tau_j\le N}
\frac{1+e^{-2\gl h(\tau_j-\tau_{j-1})-2\gl\sum_{i=(\tau_{j-1}+1)}^{\tau_j}\go_i}}2
\ind_{\{N\in\tau\}}
\right],
\end{eqnarray}
where $\ind_{\{A\}}$ is the indicator function of the event $A$.  

To see that the ``standard copolymer model'' \eqref{eq:Zsrw} is a
particular case of \eqref{eq:Z}, let $\tau:=\{n\ge0\in\N:S_n=0\}$ and
as a consequence $K(n):=\bP^{\rm SRW}(\inf\{k>0:S_k=0\}=n)$. It is
known that in this case $K(\cdot)$ satisfies \eqref{eq:K} with
$\ga=1/2$, see \cite[Ch. III]{cf:feller1} (the fact that in this case
\eqref{eq:K} holds only for $n\in2\N$, while $K(n)=0$ for $n\in 2\N+1$
due to the periodicity of the simple random walk, entails only
elementary modifications in the arguments below).  Next, observe that
if $s_i\in\{-1,+1\}$ denotes the sign of the excursion of $S$ between
the successive returns to zero $\tau_{i-1}$ and $\tau_i$, under
$\bP^{SRW}$ the sequence $\{s_i\}_{i\in\N}$ is IID and symmetric (and
independent of the sequence $\tau$).  Therefore, performing the
average on the $\{s_i\}_i$ in \eqref{eq:Zsrw} one immediately gets
\eqref{eq:Z}.

\medskip

The infinite-volume free energy is defined as
\begin{eqnarray}
  \tf(\gl,h):=\lim_{N\to\infty}\frac1N\bbE\log Z_{N,\go}\ge0,
\end{eqnarray}
where existence of the limit is a consequence of superadditivity of the
sequence
$\{\bbE\log Z_{N,\go}\}_N$ and the inequality $\tf\ge0$ is immediate from
\begin{eqnarray}
  Z_{N,\go}\ge \frac{K(N)}2,
\end{eqnarray}
which is easily seen inserting in the expectation in right-hand side of 
\eqref{eq:Z} the indicator function $\ind_{\{\tau_1=N\}}$.
One usually defines the critical line in the $(\gl,h)$ plane as
\begin{eqnarray}
  h_c(\gl):=\sup\{h\ge0:\,\tf(\gl,h)>0\}.
\end{eqnarray}
From Jensen's inequality $\bbE\log Z_{N,\go}\le \log \bbE Z_{N,\go}$
(the ``annealed bound'') one obtains the immediate inequality
\begin{eqnarray}
\label{eq:annb}
  h_c(\gl)\le \gl.
\end{eqnarray}
Indeed, one has
\begin{eqnarray}
\label{eq:Zann}
  \bbE Z_{N,\go}=
\bE\left[\prod_{j\ge 1:\tau_j\le N}
\frac{1+e^{2\gl(\gl- h)(\tau_j-\tau_{j-1})}}2
\ind_{\{N\in\tau\}}
\right],
\end{eqnarray}
from which it is not difficult to prove that
\begin{eqnarray}
\label{eq:Fann}
  \lim_{N\to\infty}\frac1N\log \bbE Z_{N,\go}=2\gl(\gl-h)\ind_{\{h<\gl\}},
\end{eqnarray}
and therefore the claim \eqref{eq:annb}.  For $h\ge\gl$,
\eqref{eq:Fann} follows from the fact that the right-hand side of
\eqref{eq:Zann} is bounded above by $1$, while the left-hand side of
\eqref{eq:Fann} is always non-negative. For $h<\gl$, just observe that
\begin{eqnarray}
\label{eq:UBann}
   \bbE Z_{N,\go}\le \bE\left[\prod_{j\ge 1:\tau_j\le N}
e^{2\gl(\gl- h)(\tau_j-\tau_{j-1})}
\ind_{\{N\in\tau\}}
\right]=\bP(N\in\tau)e^{2\gl(\gl-h)N}
\end{eqnarray}
and that
\begin{eqnarray}
  \bbE Z_{N,\go}\ge
 \bE\left[\prod_{j\ge 1:\tau_j\le N}
\frac{1+e^{2\gl(\gl- h)(\tau_j-\tau_{j-1})}}2
\ind_{\{\tau_1=N\}}
\right]=K(N)\frac{1+e^{2\gl(\gl-h)N}}2.
\end{eqnarray}
The limit \eqref{eq:Fann} is called {\sl annealed free energy}.

\subsubsection{What is known about the critical  line and its slope at the origin}

The critical point is known to satisfy the bounds
  \begin{eqnarray}
\label{eq:baundi}
    \frac{\gl}{1+\ga}\le h_c(\gl)<\gl.
  \end{eqnarray}
  The upper bound, proven recently in \cite[Th. 2.10]{cf:BGLT}, says that the
  annealed inequality \eqref{eq:annb} is strict for every $\gl$. The
  lower bound was proven in \cite{cf:BG} for the model \eqref{eq:Zsrw}
  and in the general situation \eqref{eq:Z} in \cite{cf:GB}, and is
  based on an idea by C.  Monthus \cite{cf:Monthus}.  We mention that
  (the analog of) the lower bound in \eqref{eq:baundi} was recently
  proven in \cite{cf:T_fractmom} and \cite{cf:BCT} to become optimal in
  the limit $\gl\to\infty$ for the {\sl ``reduced copolymer model''}
  introduced in \cite[Sec. 4]{cf:BG} (this is a copolymer model where the
  disorder law $\bbP$ depends on the coupling parameter $\gl$).

As already mentioned, much
  attention has been devoted to the slope of the critical curve at the
  origin, in short the ``critical slope'',
  \begin{eqnarray}
\label{eq:limite}
      \lim_{\gl\searrow0}\frac{h_c(\gl)}\gl.
  \end{eqnarray}
  Existence of such limit is not at all obvious (and indeed was proven
  \cite{cf:BdH} only in the case of the ``standard copolymer model''
  \eqref{eq:Zsrw}), but is expected to hold in general. While the
  proof in \cite{cf:BdH} was given in the case $\bbP(\go_1=\pm1)=1/2$,
  it was shown in \cite{cf:GT_ptrf} (by a much softer argument) that
  the results of \cite{cf:BdH} imply (always for the model
  \eqref{eq:Zsrw}) that the slope exists and is the same in the
  Gaussian case we are considering here. Moreover, the critical slope
  is expected to be a function only of $\ga$ and not of the full
  $K(\cdot)$, at least for $0<\ga<1$, and to be independent of the
  choice of the disorder law $\bbP$, as long as the $\go_n$'s are IID,
  with finite exponential moments, centered and of variance $1$. In
  contrast, it is known that the critical curve $\gl\mapsto h_c(\gl)$
  {\sl does} in general depend on the details of $K(\cdot)$ (this
  follows from \cite[Prop.  2.11]{cf:BGLT}) and of course on the
  disorder law $\bbP$.  The belief in the {\sl universality of the
    critical slope} is supported by the result of \cite{cf:BdH} which,
  beyond proving that the limit \eqref{eq:limite} exists, identifies
  it with the critical slope of a continuous copolymer model, where
  the simple random walk $S$ is replaced by a Brownian motion, and the
  $\go_n$'s by a white noise.

Until
  recently, nothing was known about the value of the critical slope,
  except for
\begin{eqnarray}
\label{eq:liminfsup}
  \frac1{1+\ga}\le \liminf_{\gl\searrow0}\frac{h_c(\gl)}\gl\le 
 \limsup_{\gl\searrow0}\frac{h_c(\gl)}\gl\le 1,
\end{eqnarray}
which follows from \eqref{eq:annb} and from the lower bound in
\eqref{eq:baundi} (note that the strict upper bound \eqref{eq:baundi}
does not imply a strict upper bound on the slope).
None of these bounds
is believed to be optimal. In particular, as we mentioned in the
introduction, for the standard copolymer model \eqref{eq:Zsrw}
numerical simulations \cite{cf:CGG} suggest a value around $0.83$ for
the slope.  This situation was much improved in \cite{cf:BGLT}: if
$\ga>1$, then \cite[Ths. 2.9 and 2.10]{cf:BGLT}
  \begin{eqnarray}
    \max\left(\frac1{\sqrt{1+\ga}},\frac12\right)\le\liminf_{\gl\searrow0}\frac{h_c(\gl)}\gl\le    
\limsup_{\gl\searrow0}\frac{h_c(\gl)}\gl<1.
  \end{eqnarray}
  Note that $\ga>1$ and $\ga\le 1$ are profoundly different
  situations: the inter-arrival times of the renewal process have finite mean
  in the former case and infinite mean in the latter.  Moreover, it
  was proven in \cite[Th. 2.10]{cf:BGLT} that there exists $\ga_0<1$
  (which can be estimated to be around $0.65$), such that
  if $\ga\ge\ga_0$
  \begin{eqnarray}
    \liminf_{\gl\searrow0}\frac{h_c(\gl)}\gl>\frac1{1+\ga}.
  \end{eqnarray}
Note that this does not cover the case of the standard copolymer model
\eqref{eq:Zsrw}, for which $\ga=1/2$.

\subsubsection{A new upper bound on the critical slope}
Our main result  is that the upper bound in 
\eqref{eq:liminfsup} is always strict:
\begin{theorem}
\label{th:slope}
For every $\ga>0$ there exists $\rho(\ga)\in(1/(1+\ga),1)$ such that,
whenever $K(\cdot)$ satisfies \eqref{eq:K},
\begin{eqnarray}
  \label{eq:slope}
  \limsup_{\gl\searrow0}\frac{h_c(\gl)}\gl\le 
 \rho(\ga).
\end{eqnarray}
\end{theorem}
It is interesting to note that the upper bound \eqref{eq:slope} depends
only on the exponent $\ga$ and not the details of $K(\cdot)$.
This is coherent with the mentioned belief in universality of the 
slope. 

The new idea which allows to go beyond the results of \cite[Th. 2.10]{cf:BGLT} is 
to bound above the fractional moments of $Z_{N,\go}$ in two steps:
\begin{enumerate}
\item
first we chop the system 
into blocks of size $k$, the correlation length of the annealed model, and we decompose
$Z_{N,\go}$ according to which of the blocks contain points of $\tau$ 
\item
 only at that
point we apply the inequality \eqref{eq:ineqgamma}, where each of the $a_i$ corresponds to one of
the pieces into which the partition function has been decomposed.
\end{enumerate}

\begin{rem}\rm 
  Theorem \ref{th:slope} holds in the more general situation where
  $\go$ is a sequence of IID random variables with finite exponential
  moments and normalized so that $\bbE\,\go_1=0,\bbE\, \go_1^2=1$. We
  state the result and give the proof only in the Gaussian case simply
  to keep technicalities at a minimum. The extension to the general
  disorder law can be obtained following the lines of \cite[Sec.
  4.4]{cf:BGLT}.
\end{rem}
\section{Proof of Theorem \ref{th:slope}.}

Fix  $\ga>0$, $1/(1+\ga)<\gga<1$ and define
\begin{eqnarray}
  \label{eq:A}
  A_{N,\gamma}:=\bbE(Z_{N,\go}^\gga).
\end{eqnarray}
The reason why we restrict to $\gga>1/(1+\ga)$ will be clear after
\eqref{eq:condgep}.
From now on we take $h=\rho \gl$, where
the value of 
\begin{eqnarray}
\label{eq:rho<1}
  \rho:=\rho(
\ga,\gga)\in(\gamma,1) 
\end{eqnarray}
will be chosen close
to $1$ later. Let
\begin{eqnarray}
  k:=\left\lfloor \frac1{\gl^2(1-\rho)}\right\rfloor
\end{eqnarray}
and note that, irrespective of how $\rho$ is chosen, $k$ can be made
arbitrarily large choosing $\gl$ small (which is no restriction since
in Theorem \ref{th:slope} we are interested in the limit $\gl\searrow
0$).  One sees from \eqref{eq:Fann} that, apart from an inessential
factor $2$, $k$ is just the inverse of the annealed free energy, i.e.,
\begin{eqnarray}
  k=\left\lfloor\frac2{\lim_{N\to\infty}(1/N)\log \bbE Z_{N,\go}}\right\rfloor.
\end{eqnarray}

We will show that, if $\rho(\ga,\gga)$ is sufficiently close to
$1$, there exists $\gl_0:=\gl_0(\gga,K(\cdot))>0$ such that for
$0<\gl<\gl_0$ there exists $c:=c(\gga,\gl,K(\cdot))<\infty$ such that
\begin{eqnarray}
\label{eq:auxA}
  A_{N,\gamma}\le c\,\left[K(N/k)\right]^\gamma
\end{eqnarray}
for every $N\in k\N$. In particular, by Jensen's inequality
and the fact that the sequence $\{(1/N)\bbE\log
Z_{N,\go}\}_N$ has a non-negative limit,
\begin{eqnarray}
  \tf(\gl,\rho(\ga,\gga)\gl)=\lim_{n\to\infty}
\frac1{nk}\bbE\log Z_{nk,\go}\le \lim_{n\to\infty}
\frac1{nk\gamma}\log A_{nk,\gamma}=0
\end{eqnarray}
for $\gl< \gl_0$.
 This implies \eqref{eq:slope} with 
\begin{eqnarray}
   \rho(\ga):=\inf_{1/(1+\ga)<\gga<1}\rho(\ga,\gga).
\end{eqnarray}

\begin{figure}[h]
\begin{center}
\leavevmode
\epsfysize =4 cm
\psfragscanon
\psfrag{0}[c]{$0=i_0=j_0$}
\psfrag{B1}[c]{$B_1$}
\psfrag{B2}[c]{$B_2$}
\psfrag{B3}[c]{$B_3$}
\psfrag{B4}[c]{$B_4$}
\psfrag{Z1}[c]{$Z_{n_1,j_1}$}
\psfrag{Z2}[c]{$Z_{n_\ell,N}$}
\psfrag{n1}[c]{ $n_1$}
\psfrag{j1}[c]{ $j_1$}
\psfrag{nl}[c]{$n_\ell$}
\psfrag{N}[c]{$N$}
\psfrag{t2}[c]{$\tau_2$}
\psfrag{t3}[c]{$\tau_3$}
\psfrag{t4}[c]{$\tau_4$}
\epsfbox{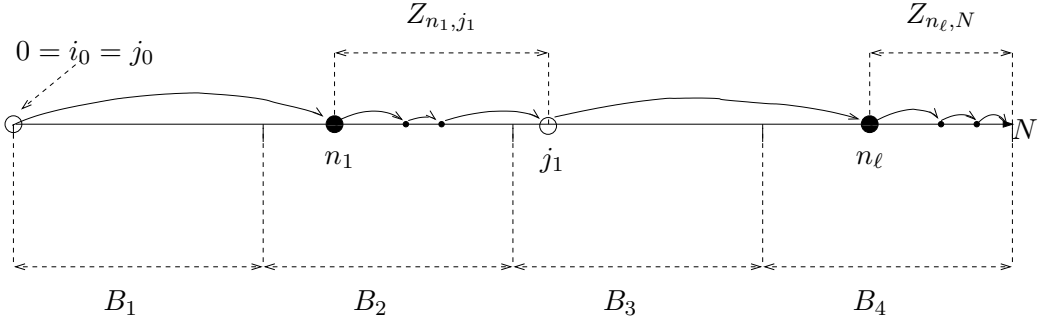}
\end{center}
\caption{\label{fig:decompos} A typical configuration which
  contributes to $\hat Z_\go^{(i_1,\ldots,i_\ell)}$ defined in
  \eqref{eq:Zhat}.  In this example we have a number $N/k=4$ of
  blocks, $\ell=2$ and $i_1=2$, while $i_\ell=N/k=4$ by definition
  (cf. \eqref{eq:dec1}).  Big black dots denote the $n_i$'s, white
  dots denote the $j_i$'s, while small dots are all the other points
  of $\tau$.  Note that $j_1-n_1<k$, as it should, and that there is
  no point of $\tau$ between a white dot and the next big black dot.
  In this example, the set $M$ of \eqref{eq:M} is $\{2,3,4\}$, and as
  a consequence $W$ defined in \eqref{eq:W} is $B_2\cup B_3\cup B_4$.
  In words: $n_1$ is the first point of $\tau$ after $0$, $j_1$ is the
  last point of $\tau$ which does not exceed $n_1+k-1$, $n_2$ is the
  first point after $j_1$, and so on.  The index of the block
  containing $n_r$ defines $i_r$. }
\end{figure}

From now on we  assume that $(N/k)$ is integer and we divide the interval
$\{1,\ldots,N\}$ into blocks
\begin{eqnarray}
  \label{eq:blocks}
  B_i:=\{(i-1)k+1,(i-1)k+2,\ldots,ik\}\;\mbox{with\;\;}i=1,\ldots,N/k.
\end{eqnarray}
Set $Z_{i,j}:=Z_{(j-i),\theta^i\go}$ (with the convention
$Z_{i,i}=1$), where $\theta$ is the left shift operator:
$(\theta^a\go)_i:=\go_{i+a}$ for $i,a\in\N$.  We have then the
identity (see Fig. \ref{fig:decompos})

\begin{eqnarray}
\label{eq:dec1}
  Z_{N,\go}=\sum_{\ell=1}^{N/k}\sum_{i_0:=0<i_1<\ldots<i_\ell=N/k}
\hat Z_\go^{(i_1,\ldots,i_\ell)}
\end{eqnarray}
where
\begin{eqnarray}
  \label{eq:Zhat}
  \hat Z_\go^{(i_1,\ldots,i_\ell)}&:=&
  \sum_{n_1\in B_{i_1}}\sum_{j_1=n_1}^{n_1+k-1}\sumtwo{n_2\in B_{i_2}:}{n_2\ge
    n_1+k}\sum_{j_2=n_2}^{n_2+k-1}\ldots \sumtwo{n_{\ell-1}\in B_{i_{\ell-1}}:}
{n_{\ell-1}\ge
    n_{\ell-2}+k}\sum_{j_{\ell-1}=n_{\ell-1}}^{n_{\ell-1}+k-1}\sumtwo
{n_\ell\in B_{N/k}:}{n_\ell\ge n_{\ell-1}+k}\\\nonumber
&&\times \varphi\left((0,n_1]\right)K(n_1)Z_{n_1,j_1}\varphi((j_1,n_2])K(n_2-j_1)
Z_{n_2,j_2}\times\ldots \\\nonumber
&&\times \varphi((j_{\ell-1},n_\ell])K(n_\ell-j_{\ell-1})Z_{n_\ell,N}
\end{eqnarray}
and, for $I\subset \N$,
\begin{eqnarray}
  \label{eq:phi}
  \varphi(I):=\frac{1+e^{-2\gl h |I|-2\gl\sum_{i\in I}\go_i}}2.
\end{eqnarray}
We have then, using the inequality
\begin{eqnarray}
  \label{eq:ineqgamma}
\left(  \sum_{i=1}^n a_i\right)^\gga\le \sum_{i=1}^n a_i^\gga
\end{eqnarray}
which holds for $0\le\gga\le1$ and $a_i\ge0$,
\begin{eqnarray}
  A_{N,\gamma}\le \sum_{\ell=1}^{N/k}\sum_{i_0:=0<i_1<\ldots<i_\ell=N/k}
\bbE\left[(\hat Z_\go^{(i_1,\ldots,i_\ell)})^\gamma\right].
\end{eqnarray}
Define
\begin{eqnarray}
  \label{eq:M}
  M:=M(i_1,\ldots,i_\ell):=\{i_1,i_2,\ldots,i_\ell\}\cup \{i_1+1,i_2+1,\ldots,
i_{\ell-1}+1\},
\end{eqnarray}
\begin{eqnarray}
  \label{eq:W}
  W:=W(i_1,\ldots,i_\ell):=\cup_{u\in M}B_u
\end{eqnarray}
and note that $1\le |M|< 2\ell$. With the conventions of Fig.
\ref{fig:decompos}, $W$ is the union of the blocks $B_i$ which either
contain a big black dot or such that $B_{i-1}$ contains a big black dot.
Note also that, for every $r$, the interval $[n_r,j_r]$ is a subset of $W$.

We want first of all to show that the $\varphi$'s can be effectively
replaced by constants.  To this purpose, we use the inequality
\begin{eqnarray}
  \label{eq:phiphi}
  \varphi((j_r,n_{r+1}])\le 2\varphi((j_r,n_{r+1}]\setminus W)
\,\varphi\left((j_r,n_{r+1}]\cap W\right)
\end{eqnarray}
(with the convention that $\varphi(\emptyset)=1$ and $j_0:=0$),
where $W$ was defined in \eqref{eq:W}.
This is simply due to the fact that, if $I_1$ and $I_2$ are two disjoint
subsets of $\N$, one has $\varphi(I_1\cup I_2)\le 2\varphi(I_1)\varphi(I_2)$.

We note that the two factors in the right-hand side of \eqref{eq:phiphi} are
independent random variables. Moreover, since
\begin{eqnarray}
  (j_r,n_{r+1}]\setminus W=\cup_{i:i_r+1<i<i_{r+1}} B_i,
\end{eqnarray}
we observe that the law of
\begin{eqnarray}
\label{eq:prod}
  \prod_{r=0}^{\ell-1}
\varphi\left((j_r,n_{r+1}]\setminus W\right)
\end{eqnarray}
depends only on $(i_1,\ldots,i_\ell)$ and not on the $n_r$'s and
$j_r$'s, and that, once $(i_1,\ldots,i_\ell)$ is fixed,
\eqref{eq:prod} is the product of $\ell$ independent random variables.
As a consequence,
\begin{eqnarray}
  A_{N,\gamma}\le \sum_{\ell=1}^{N/k}\sum_{i_0:=0<i_1<\ldots<i_\ell=N/k} 2^{\gamma\ell}
\prod_{r=0}^{\ell-1}\bbE \left(
\varphi\left((j_r,n_{r+1}]\setminus W\right)^\gamma\right)
\bbE \left[\left(\check Z^{(i_1,\ldots,i_\ell)}_\go\right)^\gamma\right],
\end{eqnarray}
where
\begin{eqnarray}
  \label{eq:Zcheck}
  \check Z^{(i_1,\ldots,i_\ell)}_\go&:=&
 \sum_{n_1\in B_{i_1}}\sum_{j_1=n_1}^{n_1+k-1}\sumtwo{n_2\in B_{i_2}:}{n_2\ge
    n_1+k}\sum_{j_2=n_2}^{n_2+k-1}\ldots \sumtwo{n_{\ell-1}\in B_{i_{\ell-1}}:}
{n_{\ell-1}\ge
    n_{\ell-2}+k}\sum_{j_{\ell-1}=n_{\ell-1}}^{n_{\ell-1}+k-1}\sumtwo
{n_\ell\in B_{N/k}:}{n_\ell\ge n_{\ell-1}+k}\\\nonumber
&&\times\varphi\left((0,n_1]\cap W\right)K(n_1)Z_{n_1,j_1}
\varphi((j_1,n_2]\cap W)K(n_2-j_1)
Z_{n_2,j_2}\times\ldots\\\nonumber 
&&
\times 
\varphi((j_{\ell-1},n_\ell]\cap W)K(n_\ell-j_{\ell-1})Z_{n_\ell,N}.
\end{eqnarray}
Thanks to \eqref{eq:ineqgamma}
and to the choice $h=\rho\gl$, for
every $I\subset \N$
\begin{eqnarray}
  \bbE (\varphi(I)^\gamma)\le \frac{1+e^{2\gl^2\gamma(\gamma-\rho)|I|}}
{2^\gamma}\le 2^{1-\gamma},
\end{eqnarray}
where the second inequality is implied by our assumption $\gamma<\rho$,
cf. \eqref{eq:rho<1}.
Then,
\begin{eqnarray}
  A_{N,\gamma}\le \sum_{\ell=1}^{N/k}\sum_{i_0:=0<i_1<\ldots<i_\ell=N/k} 2^\ell
\bbE \left[\left(\check Z^{(i_1,\ldots,i_\ell)}_\go\right)^\gamma\right].
\end{eqnarray}
In order to estimate
the remaining average, we use H\"older's inequality with $p=1/\gamma$
and $q=1/(1-\gamma)$:
\begin{eqnarray}
  \label{eq:Hold}
  \bbE \left[\left(\check Z^{(i_1,\ldots,i_\ell)}_\go\right)^\gamma\right]
=\tilde \bbE\left[\left(\check Z^{(i_1,\ldots,i_\ell)}_\go\right)^\gamma
\frac{\dd \bbP}{\dd \tilde \bbP}(\go)\right]
\le \left(
\tilde \bbE \check Z^{(i_1,\ldots,i_\ell)}_\go\right)^\gamma
\left(\tilde \bbE \left(\frac{\dd \bbP}{\dd \tilde \bbP}\right)^{1/(1-\gamma)}
\right)^{1-\gamma}
\end{eqnarray}
where, under the modified law $\tilde \bbE:=\tilde
\bbE^{(i_1,\ldots,i_\ell)}$, the $\{\go_i\}_{i\in\N}$ are still
Gaussian, independent and of variance $1$, but $\go_i$ has average
$1/\sqrt k$ if $i\in W(i_1,\ldots,i_\ell)$, while $\go_i$ has average
$0$, as under $\bbE$, if $i\notin W(i_1,\ldots,i_\ell)$. Since
$\tilde\bbE$ is still a product measure, it is immediate to check that
\begin{eqnarray}
  \left(\tilde \bbE \left(\frac{\dd \bbP}{\dd \tilde \bbP}\right)^{1/(1-\gamma)}
\right)^{1-\gamma}
=\left[\bbE\, e^{\left(\go_1/\sqrt k+1/(2k)\right)
\frac\gga{1-\gga}}
\right]^{(1-\gga)|W|}=e^{\frac{\gamma |W|}{2k(1-\gga)}}
\le e^{\frac\gamma{1-\gamma}\ell},
\end{eqnarray}
where we used the fact that $|W|=k|M|$ and $|M|\le 2\ell$.  Next, we observe that
\begin{eqnarray}
\label{eq:mnk}
  \tilde \bbE\varphi((j_r,n_{r+1}]\cap W )=
\frac{1+\tilde\bbE e^{-2\gl h |I|-2\gl \sum_{i\in I}\go_i}}2
\end{eqnarray}
with $I=(j_r,n_{r+1}]\cap W$. Thanks to the definition of 
$k$ and to $h=\rho \gl$, \eqref{eq:mnk} equals
\begin{eqnarray}
  \frac{1+e^{2\gl^2|I|(1-\rho)-2\gl|I|/\sqrt k}}2\le
\frac{1+e^{2\gl^2|I|(1-\rho)-2|I|\gl^2\sqrt{1-\rho}}}2\le 1.
\end{eqnarray}
In conclusion, we proved
\begin{eqnarray}
  A_{N,\gamma}&\le&\sum_{\ell=1}^{N/k}\sum_{i_0:=0<i_1<\ldots<i_\ell=N/k}
\left(2\, e^{\frac\gamma{1-\gamma}}\right)^\ell\\
\label{eq:sqbr}
&\times&
\left[\sum_{n_1\in B_{i_1}}\ldots \sumtwo{n_\ell\in B_{N/k}:}{n_\ell\ge 
n_{\ell-1}+k} K(n_1)\ldots K(n_\ell-j_{\ell-1})
U(j_1-n_1)\ldots U(N-n_\ell)\right]^\gamma
\end{eqnarray}
where, since $[n_r,j_r]$ is a subset of $W$ as observed above,
\begin{eqnarray}
  \label{eq:U}
  U(n):= \bbE Z_{n,\go}(\gl,\rho\gl+1/\sqrt k)
\end{eqnarray}
with the convention that $U(0):=1$.  In \eqref{eq:sqbr} we used
independence of $Z_{n_r,j_r}$ for different $r$'s (recall that
$\tilde\bbE$ is a product measure) to factorize the expectation.  The heart
of the proof of Theorem \ref{th:slope} is the following:
\begin{lemma}
\label{th:lemmagep}
There exists $\gl_0(\gga,K(\cdot))>0$ such that the following holds
for $\gl<\gl_0$.  If, for some $\gep>0$,
  \begin{eqnarray}
    \label{eq:cond1}
    \sum_{j=0}^{k-1}U(j)\le \gep\, Y(k,{\mathcal C}_K,\ga):=\gep\times \left\{
      \begin{array}{lll}
        \frac{k^\ga}{{{\mathcal C}_K}}  & {\rm if} & \ga<1 \\
 & & \\
\frac{2k}{{{\mathcal C}_K} \log k} & {\rm if} & \ga=1\\
 & & \\
 k & {\rm if} & \ga>1
      \end{array}
\right.
  \end{eqnarray}
and
  \begin{eqnarray}
    \label{eq:cond2}
    \sum_{j=0}^{k-1}\sum_{n\ge k}U(j)K(n-j)\le \gep
  \end{eqnarray}
then  the quantity in square brackets in \eqref{eq:sqbr} is bounded above
by
\begin{eqnarray}
  \label{eq:concllemma}
  C_1\, \gep^\ell\,3^{(3+2\ga)\ell} \prod_{r=1}^\ell \frac1{(i_r-i_{r-1})^{1+\ga}}
\end{eqnarray}
where $C_1:=C_1(\gep, k,K(\cdot))<\infty$.
\end{lemma}
Here and in the following, the positive and finite constants $C_i, i\ge1$
depend only on the arguments which are explicitly indicated, while
${{\mathcal C}_K}$ is the same constant which appears in \eqref{eq:K}.  \medskip

Assume that Lemma \ref{th:lemmagep} is true and that
\eqref{eq:cond1}-\eqref{eq:cond2} are satisfied. Then,
\begin{eqnarray}
\label{eq:330}
A_{N,\gamma}\le C_1(\gep, k,K(\cdot))^\gga
\sum_{\ell=1}^{N/k}\sum_{i_0:=0<i_1<\ldots<i_\ell=N/k}&&
\left(3^{\gamma(3+2\ga)}\,2\,e^{\frac\gamma{1-\gamma}}\right)^\ell \gep^{\gga \ell}
\\\nonumber
&&\times\prod_{r=1}^\ell \frac1{(i_r-i_{r-1})^{(1+\ga)\gga}}.
\end{eqnarray}
If moreover  $\gep$ satisfies
\begin{eqnarray}
  \label{eq:condgep}
\xi:=3^{\gga(3+2\ga)}\,2\,e^{\frac\gamma{1-\gamma}}\,\gep^{\gamma} \sum_{n\in\N}\frac1{n^{(1+\ga)\gga}}< 1
\end{eqnarray}
then it follows from \cite[Th. A.4]{cf:GB} that
\begin{eqnarray}
\label{eq:335}
  A_{N,\gamma}\le C_2(\gep,\gga,k,K(\cdot)) {(N/k)^{-(1+\ga)\gga}}
\end{eqnarray}
for every $N\in k\N$. Indeed, the sum in the right-hand
side of \eqref{eq:330} is nothing but the partition function of a
homogeneous pinning model \cite[Ch. 2]{cf:GB} of length $N/k$
 with pinning parameter $\xi$
such that the system is in the delocalized phase
(this is encoded in \eqref{eq:condgep}).
More precisely: to obtain \eqref{eq:335} it is sufficient to apply
Proposition \ref{th:A4} below, 
with $\ga$ replaced by $(1+\ga)\gga-1>0$ and $K(\cdot)$ replaced
by
\begin{eqnarray}
  \hat K(n)=\frac 1{n^{(1+\ga)\gga}}\left(\sum_{j\in\N}\frac 1{j^{(1+\ga)\gga}}
\right)^{-1}.
\end{eqnarray}
\begin{proposition}\cite[Th. A.4]{cf:GB}
\label{th:A4}
  If $K(\cdot)$ is a probability on $\N$ which satisfies \eqref{eq:K}
  for some $\alpha>0$, then for every $\xi<1$ there exists 
$c=c(K(\cdot),\xi)$ such that for every $N\in\N$
\begin{eqnarray}
  \sum_{\ell=1}^N\sum_{0=i_0<\ldots<i_\ell=N}\xi^\ell \prod_{r=1}^\ell
K(i_r-i_{r-1})\le c\,K(N).
\end{eqnarray}
\end{proposition}

Inequality \eqref{eq:auxA} is then proven; note that $C_2$ depends on
$\gl$ through $k$.  The condition $\gga>1/(1+\ga)$ which we required
since the beginning guarantees that the sum in \eqref{eq:condgep}
converges.  Note that \eqref{eq:condgep} depends on $K(\cdot)$ only
through $\ga$; this is important since we want $\rho$ in
\eqref{eq:slope} to depend only on $\ga$ and not on the whole
$K(\cdot)$.

\bigskip

To conclude the proof of Theorem \ref{th:slope}, we still have to
prove Lemma \ref{th:lemmagep} and to show that
\eqref{eq:cond1}-\eqref{eq:cond2} can be satisfied, with $\gep:=
\gep(\ga,\gga)$ satisfying \eqref{eq:condgep}, for all $\gl\le \gl_0$ and 
$h= \gl \rho(\ga,\gga)$,
 if $\gl_0(\gga, K(\cdot))$ is sufficiently small and $\rho(\ga,\gga)$ is close to $1$
 (see also Remark \ref{rem:flow} below).

\medskip

{\sl Proof of Lemma \ref{th:lemmagep}.}
First of all, we get rid of $U(N-n_\ell)$ and effectively we replace 
$n_\ell$ by $N$: the quantity in square brackets in \eqref{eq:sqbr}
is upper bounded by
\begin{eqnarray}
\label{eq:t1}
C_3(k,K(\cdot))
\sum_{n_1\in B_{i_1}}\ldots \sum_{j_{\ell-1}=n_{\ell-1}}^{n_{\ell-1}+k-1}&& 
K(n_1)\ldots K(n_{\ell-1}-j_{\ell-2})K(N-j_{\ell-1})\\\nonumber
&& \times U(j_1-n_1)\ldots U(j_{\ell-1}-n_{\ell-1}).
\end{eqnarray}
Explicitly, one may take
\begin{eqnarray}
\label{eq:C3}
  C_3(k,K(\cdot))=k \;U(N-n_\ell)\supthree{j,n,N:}{0<j<n}{N-k\le n\le N} \frac{K(n-j)}{K(N-j)}:
\end{eqnarray}
one has
$$
U(N-n_\ell)\le \bbE Z_{N-n_\ell,\go}(\gl,\rho\gl)\le e^{2\gl^2 k}
$$
(recall \eqref{eq:UBann} and $N-n_\ell\le k$), while the supremum in
\eqref{eq:C3} is easily seen from \eqref{eq:K} to depend only on $k$
and $K(\cdot)$.

Recall that by convention $i_0=j_0=0$, and let also from now on 
$n_\ell:=N$ and $n_0:=0$. We do the following:
\begin{itemize}
\item for every $1\le r\le \ell$ such that $i_r>i_{r-1}+2$
(which guarantees that between $j_{r-1}$ and $n_r$ there is at least
one full block), we use
  \begin{eqnarray}
    K(n_r-j_{r-1})\le 3^{\ga+2} \frac {{\mathcal C}_K}{(i_r-i_{r-1})^{1+\ga}k^{1+\ga}}.
  \end{eqnarray}
This is true under the assumption that $k\ge k_0$
with $k_0(\gga,K(\cdot))$ large enough, i.e., $\gl\le
\gl_0(\gga,K(\cdot))$ with $\gl_0$ small, since
$n_r-j_{r-1}\ge k(i_r-i_{r-1}-2)$ and 
$\sup_{n>2}(n/(n-2))^{1+\ga}=3^{1+\ga}$.
\item for every $1\le r\le \ell$ such that $i_r\le i_{r-1}+2$, we
  leave $K(n_r-j_{r-1})$ as it is.
\end{itemize}
Then, \eqref{eq:t1} is bounded above by 
\begin{eqnarray}
  \label{eq:t2}
 && C_3(k,K(\cdot))\frac{3^{(\ga+2)\ell}}{k^{(1+\ga)|J|}}\prod_{r\in J}
\frac {{\mathcal C}_K}{(i_r-i_{r-1})^{(1+\ga)}}\\\nonumber
&& \times\sum_{n_1\in B_{i_1}}\ldots \sum_{j_{\ell-1}=n_{\ell-1}}^{n_{\ell-1}+k-1}
\left(\prod_{r\in\{1,\ldots,\ell\}\setminus J}
K(n_r-j_{r-1})\right)U(j_1-n_1)\ldots U(j_{\ell-1}-n_{\ell-1}),
\end{eqnarray}
where 
\begin{eqnarray}
  \label{eq:J}
  J:=J(i_1,\ldots,i_\ell):=\{1\le r\le \ell:i_r>i_{r-1}+2\}.
\end{eqnarray}
Now we can sum over $j_r,n_r,$ $1\le r<\ell$, using the two assumptions
\eqref{eq:cond1}-\eqref{eq:cond2}. We do this in three steps:
\begin{itemize}
\item First, for every $r\in J$ we sum over the allowed values of
  $j_{r-1}$ (provided that $r>1$, otherwise there is no
  $j_{r-1}$ to sum over) using \eqref{eq:cond1} and the constraint $0\le
  j_{r-1}-n_{r-1}<k$.  The sum over all such $j_{r-1}$ gives at
  most
\begin{eqnarray}
C_{4}(\gep,k,K(\cdot))
\gep^{|J|}Y(k,{\mathcal C}_K,\ga)^{|J|}
\end{eqnarray}
where $C_{4}$ takes care of the fact that possibly $1\in J$.
At this point we are left with the task of summing
\begin{eqnarray}
  \prod_{r\in\{1,\ldots,\ell\}\setminus J}\left[
K(n_r-j_{r-1})U(j_{r-1}-n_{r-1})\right]
\end{eqnarray}
(observe that $U(j_0-n_0)=1$) over all the allowed values 
of $n_r,1\le r<\ell$ and of $j_{r-1}, r\in \{1,\ldots,\ell\}\setminus J$.
\item Secondly, using \eqref{eq:cond2}, we see that if
  $r\in\{1,\ldots,\ell\}\setminus J$ then the sum of
  $K(n_r-j_{r-1})U(j_{r-1}-n_{r-1})$ over the allowed values of
  $n_r$ and $j_{r-1}$ gives at most $ \gep$ (or $1$ if $r=1$).  The
  contribution from all the $n_r,j_{r-1}$ with
  $r\in\{1,\ldots,\ell\}\setminus J$ is therefore at most
\begin{eqnarray}
  C_{5}(\gep,k,K(\cdot))\gep^{\ell-|J|}.
\end{eqnarray}
This is best seen if one starts to sum on  $n_r,j_{r-1}$ for the largest 
value of  $r\in\{1,\ldots,\ell\}\setminus J$ 
and then proceeds to the second largest, and so on.
\item Finally, the sum over all the $n_r$'s with $r\in J$ is trivial
(the summand does not depend on the $n_r$'s) and gives at
  most $k^{|J|}$.
\end{itemize}
In conclusion, we have upper bounded \eqref{eq:t2} by
\begin{eqnarray}
\label{eq:lastfact}
&& C_{6}(\gep,k,K(\cdot))3^{(3+2\ga)\ell} \gep^{\ell}\prod_{r=1
}^\ell\frac1{(i_r-i_{r-1})^{1+\ga}}\left[\frac{{\mathcal C}_K Y(k,{\mathcal C}_K,\ga)}
{k^{\ga}}\right]^{|J|}.
\end{eqnarray}
If $0<\ga<1$ it is clear from the definition of $Y(k,{\mathcal C}_K,\ga)$ that
the last factor equals $1$ and \eqref{eq:concllemma}
is proven. For $\ga\ge1$, 
$$
\frac{{\mathcal C}_K Y(k,{\mathcal C}_K,\ga)}{k^\ga}
$$
can be made as small as wished with $k$ large (i.e. choosing $\gl_0$
small), so that we can again assume that the last factor in
\eqref{eq:lastfact} does not exceed $1$. Lemma \ref{th:lemmagep} is proven.
\qed


\medskip

Finally, we have: 
\begin{proposition} 
\label{th:cond12}
Let $ \gl_0(\gga,K(\cdot))>0$ be sufficiently small.
There exists  $\gep:=\gep(\ga,\gga)>0$ satisfying
  \eqref{eq:condgep} and $\rho:=\rho(\ga,\gga)\in (\gga,1)$ such that   
conditions
  \eqref{eq:cond1}-\eqref{eq:cond2} are satisfied for all $0<\gl<\gl_0$
  and $h=\rho\gl$.
\end{proposition}

{\sl Proof of Proposition \ref{th:cond12}.}  
We have by direct computation
\begin{eqnarray}
\label{eq:Uj}
  U(j)&=&\bE\left[\prod_{1\le n:\tau_n\le j}
\frac{1+e^{(2\gl^2(1-\rho)-2\gl/\sqrt k)(\tau_n-\tau_{n-1})}}2
\ind_{\{j\in\tau\}}\right]\\\nonumber
&\le& 
\bE\left[\prod_{1\le n:\tau_n\le j}
\frac{1+e^{-\frac1{k\sqrt{1-\rho}}(\tau_n-\tau_{n-1})}}2
\ind_{\{j\in\tau\}}\right]\\
\label{eq:condiz}
&=&\bP(j\in\tau)\bE\left[\left.\prod_{1\le n:\tau_n\le j}
\frac{1+e^{-\frac1{k\sqrt{1-\rho}}(\tau_n-\tau_{n-1})}}2
\right|j\in\tau\right]
\end{eqnarray}
where in the inequality we assume that 
$2\sqrt{1-\rho}<1$
(it is important that this condition does  not depend on $\gl$).
Of course, from \eqref{eq:Uj} we see that for every $j$
\begin{eqnarray}
U(j)\le \bP(j\in\tau). 
\end{eqnarray}
Moreover, we know from 
\cite[Th. B]{cf:Doney} that, if $0<\ga<1$,
\begin{eqnarray}
  \label{eq:Doney}
  \bP(j\in\tau)\stackrel{j\to\infty}\sim \frac{\ga\sin(\pi\ga)}
{\pi}\frac1{{{\mathcal C}_K} j^{1-\ga}},
\end{eqnarray}
while for $\ga=1$ 
\begin{eqnarray}
  \label{eq:a=1}
   \bP(j\in\tau)\stackrel{j\to\infty}\sim\frac{1}{{{\mathcal C}_K}\log j} 
\end{eqnarray}
(cf. for instance  \cite[Th. A.6]{cf:GB}).
For every $a\in(0,1)$  one has then for $k$ sufficiently large, i.e., for
$\gl_0$ small,
\begin{eqnarray}
  \label{eq:somma1}
  \sum_{j=0}^{ak}U(j)\le   \sum_{j=0}^{ak}\bP(j\in\tau)\le 
Y(ak,{\mathcal C}_K,\ga)
\end{eqnarray}
(we recall that $Y(\cdot,\cdot,\cdot)$ was defined in \eqref{eq:cond1}).
We need also the following fact:
\begin{lemma}
\label{th:lemma43}
For every $\ga>0$ there exists $C_7(\ga)<\infty$ such that the following holds.
If $0<\ga<1$ then, say, for $q>2$
  \begin{eqnarray}
\label{eq:fa}
    \limsup_{N\to\infty}\bE
\left[\left.
\prod_{1\le j:\tau_j\le N}\frac{1+e^{-(q/N)(\tau_j-\tau_{j-1})}}2
\right|N\in\tau
\right]
\le C_7(\ga)\frac{(\log q)^2}{q^\ga}.
  \end{eqnarray}
If  $\ga\ge1$ then, for every $q>0$,
  \begin{eqnarray}
\label{eq:faga1}
    \limsup_{N\to\infty}\bE
\left[\left.
\prod_{1\le j:\tau_j\le N}\frac{1+e^{-(q/N)(\tau_j-\tau_{j-1})}}2
\right|N\in\tau
\right]
\le C_7(\ga) e^{-q/2}.
  \end{eqnarray}
\end{lemma}
Lemma \ref{th:lemma43} was proven in \cite[Lemma 4.3]{cf:BGLT}.  We add a few
side remarks about its proof in Appendix \ref{sec:lemma43}.

The upper bound \eqref{eq:fa} is certainly not optimal, but 
it gives us an estimate which vanishes for $q\to\infty$ and
which depends only on $\ga$ and $q$, which is all we need in the following.
Let us mention that in the case of the standard copolymer model,
using the property that for every $N\in2\N$ and every $k\in\{0,2,\ldots,N\}$
\begin{eqnarray}
  \bP^{\rm SRW}\left[\left.|\{1\le n\le N:\sign(S_n)=-1\}|=k
\right|N\in\tau\right]=\frac1{(N/2)+1}
\end{eqnarray}
(see \cite[Ch. III.9]{cf:feller1}) we obtain for every $q>0$
\begin{eqnarray}
  \lim_{N\to\infty}\bE^{\rm SRW}
\left[\left.
\prod_{1\le j:\tau_j\le N}\frac{1+e^{-(q/N)(\tau_j-\tau_{j-1})}}2
\right|N\in\tau
\right]=\frac{1-e^{-q}}q.
\end{eqnarray}

\subsubsection{Proof of \eqref{eq:cond1}}
Fix $\gep:=\gep(\ga,\gga)$ which satisfies \eqref{eq:condgep} and choose 
\begin{eqnarray}
  \label{eq:agep}
a:=a(\gep,\ga):=\left\{
\begin{array}{lll}
(\gep/2)^{1/\ga}& {\rm if} &0<\ga<1\\
& & \\
\gep/4  & {\rm if} &\ga=1\\
 & & \\
\gep/2 & {\rm if} &\ga>1.
\end{array}
\right.
\end{eqnarray}
Via \eqref{eq:somma1} one finds (for $k$ sufficiently large)
\begin{eqnarray}
\label{eq:asac2}
  \sum_{j=0}^{ak}U(j)\le \frac\gep2 Y(k,{\mathcal C}_K,\ga).
\end{eqnarray}
Next we observe that, for $ak\le j<k$, choosing $\gl_0$ small and
$\rho$ sufficiently close to $1$ (how close, depending on $K(\cdot)$
only through the exponent $\ga$) we have $U(j)\le
(\gep/2)\bP(j\in\tau)$. This just follows from Lemma \ref{th:lemma43}
above (applied with $N\sim a k$ and $q=a/\sqrt{1-\rho}$) and from
\eqref{eq:condiz}, since $\gl_0$ small implies $k$ large.  As a
consequence,
\begin{eqnarray}
\label{eq:asac}
  \sum_{j=ak}^{k-1}U(j)\le \frac\gep2 Y(k,{\mathcal C}_K,\ga)
\end{eqnarray}
and \eqref{eq:cond1} follows.

\begin{rem}
\label{rem:flow}
  \rm It is probably useful to summarize the logic of the proof of
  \eqref{eq:cond1} (similar observations hold for the proof of
  \eqref{eq:cond2} below). Given $\ga$, one first fixes $1/(1+\ga)<\gga<1$, then
  $\varepsilon(\ga, \gga)$ which satisfies \eqref{eq:condgep}, then
  $a=a(\varepsilon,\ga)$ as in \eqref{eq:agep} and
  $\rho=\rho(\ga,\gga)<1$ such that the right-hand side of Eqs.
  \eqref{eq:fa}-\eqref{eq:faga1} is smaller than $(\gep/2)$ when $q$
  is replaced by $a/\sqrt{1-\rho}$. Once all these parameters are
  fixed, one chooses $\gl_0$ sufficiently small (i.e. $k$ sufficiently
  large) so that for all $\gl<\gl_0,h=\rho\gl$ the estimates
  \eqref{eq:asac2}-\eqref{eq:asac} hold.

\end{rem}

\subsubsection{Proof of \eqref{eq:cond2}}
If we choose $b:=b(\gep,\ga)$ small, $j\le bk$ implies $k-j\ge k/2$.
Therefore, 
\begin{eqnarray}
  \sum_{j=0}^{bk}\sum_{n\ge k}U(j)K(n-j)\le
C_{8}(\ga)\frac{C_{K}}{ k^\ga} \sum_{j=0}^{bk}\bP(j\in\tau)\le
\frac\gep2
\end{eqnarray}
(if $k$ is sufficiently large and $b(\gep,\ga)$ is suitably
small). As for the rest of the sum: again from Lemma \ref{th:lemma43}
and \eqref{eq:condiz}, one has $U(j)\le (\gep/2) \bP(j\in\tau)$ for
every $bk\le j<k$.  Then,
\begin{eqnarray}
   \sum_{j=bk}^{k-1}\sum_{n\ge k}U(j)K(n-j)\le
\frac\gep2 \sum_{j=0}^{k-1}\sum_{n\ge k}\bP(j\in\tau)K(n-j)=\frac\gep2.
\end{eqnarray}
In the last equality, we used the fact that
\begin{eqnarray}
\label{eq:spiegaz}
  \sum_{j=0}^{k-1}\bP(j\in\tau)K(n-j)
\end{eqnarray}
is the $\bP$-probability that the first point of $\tau$ which does not
precede $k$ equals $n$. The sum over $n\ge k$ of \eqref{eq:spiegaz} then 
clearly equals $1$, since $\tau$ is recurrent.
\qed


\appendix

\section{Remarks on the proof of Lemma \ref{th:lemma43}.}
\label{sec:lemma43}

The proof of Lemma \ref{th:lemma43} given in \cite[Lemma 4.3]{cf:BGLT} works as follows.
Let $X_N:=\max\{n=0,\ldots,N:n\in\tau\}$, i.e., the last
point of $\tau$ up to $N$.
First of all one shows that
\begin{eqnarray}
\label{eq:nocondiz}
&&\limsup_{N\to\infty} \bE
\left[\left.
\prod_{1\le j:\tau_j\le N}\frac{1+e^{-(q/N)(\tau_j-\tau_{j-1})}}2
\right|N\in\tau
\right]\\
\label{eq:nocondiz2}
&&
\label{eq:C9}
\le C_{9}(\ga) \limsup_{N\to\infty} 
 \bE
\left[
\frac{1+e^{-(q/N)(N-X_N)}}2
\prod_{1\le j:\tau_j\le N}\frac{1+e^{-(q/N)(\tau_j-\tau_{j-1})}}2
\right].
\end{eqnarray}
We detail below the proof of this inequality in order to leave no
doubts on the fact that the constant $C_9$  depends
only on $\ga$. This was not emphasized in the proof of \cite[Lemma
4.3]{cf:BGLT} since it was not needed there.  

For $\ga\ge1$, it follows from \cite[Eqs. (4.23) and (4.49)]{cf:BGLT}
that the $\limsup$ in the right-hand side of \eqref{eq:C9} is actually
a limit, and equals $\exp(-q/2)$ (the $q/8$ which appears in \cite[Eq.
(4.49)]{cf:BGLT} can be immediately improved into $q/2$). As a side
remark, the expectation in \eqref{eq:C9}, irrespective of the value of
$\ga$ and $N$, is not smaller than $\exp(-q/2)$; this just follows
from the convexity of the exponential function:
$$
\frac{1+e^{-(q/N)x}}2\ge e^{-q/(2N)\,x}.
$$
For $0<\ga<1$, the $\limsup$ in
\eqref{eq:C9} does not exceed $C_{10}(\ga)(\log q)^2/q^\ga$, as was
proven in \cite[Eq. (4.43)]{cf:BGLT}.

Finally we prove \eqref{eq:nocondiz2}, which is quite standard.
The expectation in \eqref{eq:nocondiz} is bounded above by
\begin{eqnarray}
&& \bE
\left[\left.
\frac{1+e^{-(q/N)(N/2-X_{N/2})}}2
\prod_{1\le j:\tau_j\le N/2}\frac{1+e^{-(q/N)(\tau_j-\tau_{j-1})}}2
\right|N\in\tau
\right]\\
&&=\sum_{i=0}^{N/2} \bE
\left[\left.
\frac{1+e^{-(q/N)(N/2-i)}}2
\prod_{1\le j:\tau_j\le N/2}\frac{1+e^{-(q/N)(\tau_j-\tau_{j-1})}}2
\right|X_{N/2}=i
\right]\\&&\times\bP(X_{N/2}=i|N\in\tau)
\end{eqnarray}
and \eqref{eq:nocondiz2} follows if we can prove that
\begin{eqnarray}
  \limsup_{N\to\infty}\max_{0\le i\le N/2}\frac{\bP(X_{N/2}=i|N\in\tau)}
{\bP(X_{N/2}=i)}\le C_{11}(\ga).
\end{eqnarray}
To show this, we  use repeatedly
\eqref{eq:Doney} and \eqref{eq:K}. We start from the identity
\begin{eqnarray}
  \frac{\bP(X_{N/2}=i|N\in\tau)}
{\bP(X_{N/2}=i)}=\frac{
\sum_{j=(N/2)+1}^N K(j-i)\bP(N-j\in\tau)
}{\bP(N\in\tau)\sum_{j=(N/2)+1}^\infty K(j-i)}.
\end{eqnarray}
The denominator is lower bounded, uniformly in $0\le i\le N/2$,
by
\begin{eqnarray}
  \bP(N\in\tau)\sum_{j=(N/2)+1}^\infty K(j)\ge \frac{C_{12}(\ga)}{N},
\end{eqnarray}
where the last inequality holds for $N$ sufficiently large.
As for the numerator: always for $N$ sufficiently large,
\begin{eqnarray}
  \sum_{j=(N/2)+1}^{(3/4)N} K(j-i)\bP(N-j\in\tau)\le
C_{13}\bP(N\in\tau)\sum_{j=(N/2)+1}^{\infty}K(j-i)
\end{eqnarray}
and, uniformly on $0\le i\le N/2$,
\begin{eqnarray}
   \sum_{j=(3/4)N}^{N} K(j-i)\bP(N-j\in\tau)\le C_{14}\frac{{\mathcal C}_K}{N^{1+\ga}}
\sum_{j=0}^{N}\bP(j\in\tau)\le \frac{C_{15}(\ga)}{N}.
\end{eqnarray}
\qed

\section*{Acknowledgments}
This work was partially supported by ANR, grant POLINTBIO and grant
LHMSHE. I wish to thank the anonymous referees for the careful reading of the
manuscript and for several useful comments.

\end{document}